\documentclass[11pt]{article}
\usepackage[utf8]{inputenc}
\usepackage{graphicx} 
\usepackage{amsmath,amsfonts,amssymb,amsthm}
\usepackage[mathscr]{eucal}
\usepackage{amscd}
\usepackage{tikz}
\usepackage{tikz-cd}
\usepackage{tkz-graph}
\usepackage{multicol}
\usepackage{float}
\usepackage{authblk}

\newtheorem{theorem}{THEOREM}[section]

\theoremstyle{plain}

\newtheorem{conjecture}[theorem]{CONJECTURE}

\theoremstyle{definition}

\newtheorem{example}[theorem]{EXAMPLE}
\newtheorem{remark}[theorem]{REMARK}

\newcommand*{\dif}{\,\mathrm{d}}
\newcommand*{\mo}{\,\mathrm{mod}\,}
\newcommand*{\fr}{F}

\title{The asymptotic estimation of prime ideals in imaginary quadratic fields  and  Chebyshev’s bias }
\author{Chen Lin$^1$, Chenhao Tang$^2$, Xuejun Guo$^3$\footnote{Corresponding author.}}
\affil{{\small {$^{1,3}$School of Mathematics, Nanjing University, Nanjing 210093, China\\ $^{2}$Morningside Center of Mathematics, Academy of Mathematics and Systems Science,
Chinese Academy of Sciences; University of the Chinese Academy of Sciences,
Beijing 100190, China}\\
$^1$chen.lin@smail.nju.edu.cn,  $^2$tangchenhao25@mails.ucas.ac.cn, 
$^3$guoxj@nju.edu.cn} }
\date{}

\begin{document}
\maketitle
\begin{abstract} We study the asymptotic estimation of prime ideals that satisfy certain congruence and argument conditions in imaginary quadratic fields. We also discuss the phenomenon of Chebyshev's bias in the distribution of prime ideals among different residue classes.  
\end{abstract}

\section{Introduction}
Fermat proved that every prime congruent to $1$ modulo $4$  is a sum of two squares.
Let $p$ be an odd prime satisfying $p=x_p^2+y_p^2$, where  $x_p>y_p>0$ are integers.  
 In 1919, Hecke demonstrated that $\theta_p=\arctan \frac{y_p}{x_p}$ exhibits a uniform distribution as the prime $p$ varies.

Coleman  established the result of argument equidistribution of prime ideals in  imaginary quadratic  fields in \cite{Coleman1990TheDO} and    \cite{Coleman1992}. To state Coleman's theorem, we first recall some notation. Let $K=\mathbb{Q}(\sqrt{\Delta})$ be the imaginary quadratic number field, where $\Delta$ is the fundamental discriminant of $K$. We write $\mathfrak{f},\mathfrak{a}$ for ideals and $\mathfrak{p}$ for a prime ideal. Let $N(\mathfrak{a})$ denote the norm of $\mathfrak{a}$. Given a non-zero ideal $\mathfrak{f}$, let $g=g(\mathfrak{f})$ be the number of units $\epsilon$ such that $\epsilon\equiv 1\pmod{\mathfrak{f}}$. We write $Cl(K,\mathfrak{f})$ for the ideal class group $\mo\mathfrak{f}$ and $h_{\mathfrak{f}}$ for its cardinality. Let $C$ denote an ideal class $\mo\mathfrak{f}$ and $(\xi)$ denote the principal fractional ideal generated by some $\xi\in K$. Assume that for each $C\in Cl(K,\mathfrak{f})$, there has been fixed an ideal $\mathfrak{a}_0 \in C^{-1}$. Then for $\mathfrak{a}\in C$, we can find some $\xi_{\mathfrak{a}}\in K$ satisfying $(\xi_{\mathfrak{a}})=\mathfrak{a}\mathfrak{a}_0$ and $\xi_{\mathfrak{a}}\equiv 1 \pmod{\mathfrak{f}}$. We write $\lambda(\mathfrak{a})=\left(\frac{\xi_{\mathfrak{a}}}{|\xi_{\mathfrak{a}}|}\right)^g$.

\begin{theorem}[Coleman, Theorem 2.1 of \cite{Coleman1990TheDO}]
\label{Coleman}

Given $0\leq\varphi_1 \leq\varphi_2 \leq 2\pi,0\leq y\leq x$ and an ideal class $C\in Cl(K,\mathfrak{f}),~$we define
$$
S=S(x,y,C,\varphi_1,\varphi_2)=\{\mathfrak{a}\in C\mid x-y\leq N(\mathfrak{a})\leq x,\varphi_1 \leq arg(\lambda(\mathfrak{a}))\leq \varphi_2 \}.
$$
and
$$
T=T(x,y,C,\varphi_1,\varphi_2)=\{\mathfrak{p}\in S\mid N(\mathfrak{p})=p,~ p\text{ is prime}\}.
$$
Let $\epsilon>0$ be given. We have the asymptotic result,
$$
\sum_{\mathfrak{p}\in T(x,y,C,\varphi_1,\varphi_2)}1=\frac{(\varphi_2-\varphi_1)  y}{2\pi  h_{\mathfrak{f}} \log x }\left(1+O\left(\frac{1}{\log x}\right)\right).
$$
for $\varphi_2-\varphi_1>x^{-\frac{5}{24}+\epsilon}, y>x^{\frac{19}{24}+\epsilon}, x>x_{\epsilon}$.
\end{theorem}

A special case of Theorem \ref{Coleman} shows that the prime ideals of $\mathbb{Q}(\sqrt{-1})$, whose norms are rational primes, are argument equidistributed. Inspired by Dirichlet's Theorem on the prime distribution in arithmetic progressions, we ask the following question: among these prime ideals, are those satisfying $N(\mathfrak{p})\equiv 1\pmod{8}$ (resp. $N(\mathfrak{p})\equiv 5\pmod{8}$) argument equidistributed? The main theorem of this paper answers questions of this type and provides a generalization of Coleman's work, i.e., we present the argument equidistribution and asymptotic estimation of prime ideals that satisfy a given congruence relation in imaginary quadratic fields.

\begin{theorem}
\label{main theorem}
Using the same notation as in  Theorem \ref{Coleman}, fix two positive integers $m,M$ satisfying $\mathrm{gcd}(m,M)=1$, and define 
$$
P(x,y,C,\varphi_1,\varphi_2,m,M)=\{\mathfrak{p}\in T\mid N(\mathfrak{p})\equiv  m\pmod{M}\}.
$$
Let $\epsilon>0$ be given. We have the asymptotic result,
$$
\sum_{\mathfrak{p}\in P(x,y,C,\varphi_1,\varphi_2,m,M)}1=\frac{A(m,M,C,\mathfrak{f}) (\varphi_2-\varphi_1)  y}{2\pi  h_{\mathfrak{f}} \varphi(M) \log x }\left(1+O\left(\frac{1}{\log x}\right)\right).
$$
for $\varphi_2-\varphi_1>x^{-\frac{5}{24}+\epsilon}, y>x^{\frac{19}{24}+\epsilon}, x>x_{\epsilon}$, where $\varphi$ is Euler's totient function and
$$
A(m,M,C,\mathfrak{f})=\sum_{\phi,\psi}\phi(m)\psi(C),
$$
the sum ranging over $\phi\in\widehat{(\mathbb{Z}/M\mathbb{Z})^{\times}}, \psi\in\widehat{Cl(K,\mathfrak{f})}$ satisfying $\phi(N(\mathfrak{p}))\psi([\mathfrak{p}])=1$
for almost all prime ideals $\mathfrak{p}$ of $\mathcal{O}_K$.
\end{theorem}

It is possible that for some $m$, there does not exist a prime ideal $\mathfrak{p}$ satisfying $N(\mathfrak{p})\equiv  m\pmod{M}$. Even though such $\mathfrak{p}$ exists, it may never belong to some ideal class $C$. In such cases, we have $A(m,M,C,\mathfrak{f})=0$. But excluding these cases, the complicated coefficient $A(m,M,C,\mathfrak{f})$ has an explicit number-theoretic explanation. Let $\mathbb{Q}(\zeta_M)$ be the $M$-th cyclotomic field and $K(\mathfrak{f})$ the finite abelian extension of $K$ corresponding to $Cl(K,\mathfrak{f})$ in class field theory. Then $A(m,M,C,\mathfrak{f})$ is actually equal to $[K(\mathfrak{f})\cap\mathbb{Q}(\zeta_M):\mathbb{Q}]$, which is independent of $m$ and $C$. This will be explained in Section \ref{section 2.3}.

Many results which essentially depend on the equidistribution of prime ideals can be generalized by Theorem \ref{main theorem}. As an application of Theorem \ref{main theorem}, we generalize C. Elsholtz and G. Harman's work \cite{MR3467390} on the conjectures of T. Ordowski and Z.-W. Sun \cite{sun2017conjectures}. 

T. Ordowski and Z.-W. Sun \cite{sun2017conjectures} considered all pairs of positive integers uniquely representing a prime number by a specific quadratic form, and found that the limit of certain coordinate-wise defined functions converges. C. Elsholtz and G. Harman \cite{MR3467390} discussed these conjectures and generalized the result to cover all positive definite quadratic forms that can represent prime numbers by using Coleman's result as stated above.

By Theorem 1.2, we can extend the results of C. Elsholtz and G. Harman [4]. We focus on primes satisfying the congruence condition $p \equiv m \pmod M$. We prove in Theorem 1.3 that with this condition, the limit of certain coordinate-wise defined functions remains the same.

Let $P_{m,M}=\{p\mid p\equiv  m\pmod{M}, p \text{ is prime}\}$, where $m,M$ are coprime integers. Some definitions and notations in the following theorem will be explained later in Section \ref{pf of application}. 

\begin{theorem}
    \label{main}
    Let $Q(x,y)=ax^2+bxy+cy^2~(a,b,c\in\mathbb{Z})$ be a primitive positive definite binary quadratic form with $D=-\Delta=4ac-b^2>0$. Let $m,M$ be coprime positive integers such that $P_{m,M}\cap\{Q(x,y)\mid x,y \in\mathbb{Z}\}$ has positive density in $P_{m,M}$. For every prime $p$, define 
    $$\Omega_p=\left\{(x_p,y_p)\in \mathbb{Z}^2\mid Q(x_p,y_p)=p, x_p>y_p>0\right\}.$$
Let $k$ be a positive integer. Then we have 
$$\lim_{N\to\infty}\frac{\sum\limits_{p\leq N,p\in P_{m,M}}~\sum\limits_{(x_p,y_p)\in \Omega_p}x_p^k}{\sum\limits_{p\leq N,p\in P_{m,M}}~\sum\limits_{(x_p,y_p)\in \Omega_p}y_p^k}=\frac{\int_0^\beta{s(k,\theta)}\dif\theta}{\int_0^\beta{t(k,\theta)}\dif\theta},$$
where
$$s(k,\theta)=(\sqrt{D}\cos\theta-b\sin\theta)^k,\quad t(k,\theta)=(2a)^k \sin^k\theta,$$
and
$$\beta=\left\{
\begin{aligned}
    &\frac{1}{2}\pi &\text{if }~ 2a+b=0,\\
    \arctan&\left(\sqrt{D}/(2a+b)\right) &\text{otherwise}.
\end{aligned}
\right. $$
Here we adopt the convention that $\arctan(\cdot) \in [0, \pi)$.
\end{theorem}

Chebyshev noticed a famous phenomenon that primes congruent to 3 modulo 4 predominate over those congruent to 1 modulo 4. To be more precise, let $\pi(x;q,a)=\#\{p<x\mid p\equiv a\pmod{q} \text{ and } p \text{ is a prime number}\}$. Then, the inequality $\pi(x;4,3)>\pi(x;4,1)$ holds for most values of $x$. Littlewood \cite{littlewood1914distribution} showed that the function $\pi(x;4,3)-\pi(x;4,1)$ changes sign infinitely many times. Such a phenomenon and its generalizations are called ``Chebyshev's bias". Rubinstein and Sarnak \cite{em/1048515870} studied this phenomenon and used logarithmic density to explain this phenomenon. 

We are interested in whether the coordinate-wise defined function in our theorem has similar properties.
Using the same notation as in  Theorem \ref{main}, let 
$$\fr(N;M,m)=\frac{\sum\limits_{p\leq Pr(N),p\in P_{m,M}}~\sum\limits_{(x_p,y_p)\in \Omega_p}x_p}{\sum\limits_{p\leq Pr(N),p\in P_{m,M}}~\sum\limits_{(x_p,y_p)\in \Omega_p}y_p},$$
where $Pr(N)$ denotes the $N$-th prime number. Now using certain congruence conditions, we can divide all the primes representable by a given quadratic form into two sets, and investigate the function defined above with respect to the two sets.

For $Q(x,y)=x^2+y^2$, we compute $\fr(N;8,1)$ and $\fr(N;8,5)$ with $N$ up to $5\times 10^6$. For $Q(x,y)=x^2+xy+y^2$, we compute $\fr(N;12,1)$ and $\fr(N;12,7)$ with $N$ up to $1\times 10^6$. As shown in the figures in Subsection \ref{fr}, in each case, although the two functions converge to the same limit, they exhibit oscillations, repeatedly intersect and alternately surpass each other.

Since only oscillations and intersections matter, we consider the following new function
$$R(N;M,m)=\frac{\fr(N;M,m)}{\fr(N)},$$
where $\fr(N)=\frac{\sum_{p\leq Pr(N)}~\sum_{(x_p,y_p)\in \Omega_p}x_p}{\sum_{p\leq Pr(N)}~\sum_{(x_p,y_p)\in \Omega_p}y_p}$.
We compute $R(N;8,1)$ and $R(N;8,5)$ with $N$ up to $5\times 10^6$, and the result is shown in Figure \ref{R1-R2}. Since $\fr(N;8,1),\fr(N;8,5)$ and $\fr(N)$ have the same limit as $N\rightarrow \infty$, Figure \ref{R1-R2} shows very nice symmetry. 

The oscillations and repeated intersections in these figures suggest the Chebyshev's bias phenomenon and infinitely many sign changes of their difference (similar to Littlewood's result) in this case. On the other hand, these figures also show similarity to the phenomenon called ``murmuration". ``Murmuration" is a phenomenon involving oscillatory behavior \cite{lee2024murmurations}.

Motivated by a conjecture of Devin \cite{devin2021discrepancies}, we also study the representations $p=a^2+4b^2$ separately for primes $p\equiv1\pmod 8$ and $p\equiv5\pmod 8$. We compare the number of representations satisfying $|a|>|2b|$ with those satisfying $|a|<|2b|$. Our computations up to $1.5\times10^6$, shown in Figure \ref{D1D2}, suggest a possible bias toward $|a|<|2b|$. This leads to Conjecture \ref{conjecture_D1D2}. 

\section{Proof of Theorem \ref{main theorem}}
\label{cor of coleman}
In this section, we generalize Theorem \ref{Coleman} to Theorem \ref{main theorem}. In Section \ref{section 2.1}, we give the statement of our main results, but with rough remainder estimations. In Section \ref{section 2.2} we give the complete proof. The proof idea is to construct special Hecke L-functions to pick out the prime ideals satisfying given conditions. In Section \ref{section 2.3}, we reinterpret Theorem \ref{2.1}, Theorem \ref{2.2} from the view of $\check{\mathrm{C}}$ebotarev density theory. Via this view and combining with Theorem \ref{Coleman}, we improve the remainder estimation and obtain Theorem \ref{main theorem} in Section \ref{section 2.4}.

\subsection{Preliminaries and Statement of Results}
\label{section 2.1}
Let $K$ be a number field. For an integral ideal $\mathfrak{a}$ of $K$, we denote by $Cl(K,\mathfrak{a})$ the ideal class group mod $\mathfrak{a}$ and $h_{\mathfrak{a}}$ its cardinality, which is known to be finite. Take two positive integers $m,M$ with $\mathrm{gcd}(m,M)=1$, an integral ideal $\mathfrak{a}$ and an ideal class $C$  $\mo\mathfrak{a}$.

For any $x>0$, define two sets of prime ideals
$$
\mathscr{P}(x,C,\mathfrak{a})=\{\mathfrak{p}\in C\mid N(\mathfrak{p})\leq x\},
$$
$$
\mathscr{P}(x,m,M,C,\mathfrak{a})=\{\mathfrak{p}\in C\mid N(\mathfrak{p})\leq x, N(\mathfrak{p})\equiv m\pmod{M}\},
$$
where $N$ denotes the ideal norm.
\begin{theorem}
\label{2.1}
We have the asymptotic result
$$
\sum_{\mathfrak{p}\in\mathscr{P}(x,m,M,C,\mathfrak{a})}1=\frac{A(m,M,C,\mathfrak{a})  x}{h_{\mathfrak{a}} \varphi(M) \log x }+o\left(\frac{x}{\log x}\right)
$$
as $x\to +\infty$, where $\varphi$ is Euler's totient function and
$$
A(m,M,C,\mathfrak{a})=\sum_{\phi,\psi}\phi(m)\psi(C),
$$
the sum ranging over $\phi\in\widehat{(\mathbb{Z}/M\mathbb{Z})^{\times}}, \psi\in\widehat{Cl(K,\mathfrak{a})}$ satisfying $\phi(N(\mathfrak{p}))\psi([\mathfrak{p}])=1$
for almost all prime ideals $\mathfrak{p}$ of $\mathcal{O}_K$.
\end{theorem}

The proof of Theorem \ref{2.1} is provided in Section \ref{section 2.2} below.

We now focus on imaginary quadratic fields. Let $K=\mathbb{Q}(\sqrt{\Delta})$ be an imaginary quadratic field where $\Delta$ is the fundamental discriminant of $K$.
We denote by $w$ the number of roots of unity in $K$; $w=4$ if $\Delta=-4$, $w=6$ if $\Delta=-3$ and $w=2$ if $\Delta\neq -4,-3$.

Let $I(K)$ be the fractional ideal group of $K$ and $P(K)$ be the principal fractional ideal group. Define a homomorphism
\[
    \lambda:P(K)\longrightarrow S^{1}, (a)\longmapsto\left(\frac{a}{|a|}\right)^{w}.
    \]
To consider the argument equidistribution of prime ideals, we first extend $\lambda$ to $I(K)$. By the structure theorem of finite abelian groups, we have 
\[
    Cl(K)\cong\prod_{i=1}^{r}C_{i},
\]
where $C_{i}$ is a cyclic group with cardinality $s_i$. For each $i$, choose an integral ideal $\mathfrak{b}_i$ such that $\overline{\mathfrak{b}_{i}}$ is the generator of $C_i$. Take $\lambda(\mathfrak{b}_i)=\lambda(\mathfrak{b}_{i}^{s_i})^{\frac{1}{s_i}}$. Then for any $\mathfrak{a}\in I(K), \mathfrak{a}=(a)\prod_{i=1}^{r}(\mathfrak{b}_i)^{t_i}$, take $\lambda(\mathfrak{a})=\lambda((a))\prod_{i=1}^{r}\lambda(\mathfrak{b}_i)^{t_i}$. It is easy to check the definition is well-defined. 

Take $m,M,\mathfrak{a},C$ as at the beginning of this section. Denote by $\overline{C}$ the image of $C$ under the canonical projection $Cl(K,\mathfrak{a})\to Cl(K)$. Fix a fractional ideal $\mathfrak{a}_0\in \overline{C}^{-1}$.
Then for any $x>0$ and $ 0\leq\varphi_1 \leq\varphi_2 \leq 2\pi$,  define
\[\mathscr{P}(x,m,M,C,\mathfrak{a},\varphi_1,\varphi_2)=\{\mathfrak{p}\in \mathscr{P}(x,m,M,C,\mathfrak{a})\mid \varphi_1 \leq \mathrm{arg}\lambda(\mathfrak{p} \mathfrak{a}_0)\leq\varphi_2 \}.
\]
\begin{theorem}
\label{2.2}
    The prime ideals are equidistributed w.r.t the arguments, i.e. we have the asymptotic result
\[\sum_{\mathfrak{p}\in\mathscr{P}(x,m,M,C,\mathfrak{a},\varphi_1,\varphi_2)}1=\frac{(\varphi_2-\varphi_1)  A(m,M,C,\mathfrak{a})  x}{2\pi  h_{\mathfrak{a}} \varphi(M) \log x }+o\left(\frac{x}{\log x}\right)
\]
as $x\to +\infty$, where $\varphi$ is Euler's totient function and
$$
A(m,M,C,\mathfrak{a})=\sum_{\phi,\psi}\phi(m)\psi(C),
$$
the sum ranging over $\phi\in\widehat{(\mathbb{Z}/M\mathbb{Z})^{\times}}, \psi\in\widehat{Cl(K,\mathfrak{a})}$ satisfying $\phi(N(\mathfrak{p}))\psi([\mathfrak{p}])=1$
for almost all prime ideals $\mathfrak{p}$ of $\mathcal{O}_K$.
\end{theorem}

The proof of Theorem \ref{2.2} is also provided in Section \ref{section 2.2} below.

\begin{remark}
  If one takes the sum of the asymptotic formula over $(\mathbb{Z}/M\mathbb{Z})^{\times}$ or $Cl(K,\mathfrak{a})$ (or both) in Theorem \ref{2.1} and Theorem \ref{2.2}, one will get asymptotic results of the distribution of prime ideals with less restrictive conditions. For example, we take the sum over both $(\mathbb{Z}/M\mathbb{Z})^{\times}$ and $Cl(K,\mathfrak{a})$ in Theorem 2.1. By the orthogonality relations for characters 
    \begin{align*}
        \sum_{m\in(\mathbb{Z}/M\mathbb{Z})^{\times}}\sum_{C\in Cl(K,\mathfrak{a})}A(m,M,C,\mathfrak{a})&=\sum_{m\in(\mathbb{Z}/M\mathbb{Z})^{\times}}\sum_{C\in Cl(K,\mathfrak{a})}\sum_{\phi,\psi}\phi(m)\psi(C)\\
        &=\sum_{\phi,\psi}\sum_{m\in(\mathbb{Z}/M\mathbb{Z})^{\times}}\phi(m)\sum_{C\in Cl(K,\mathfrak{a})}\psi(C)\\
        &=\varphi(M)h_{\mathfrak{a}}.
    \end{align*}
    Thus,
    $$
    \sum_{N(\mathfrak{p})\le x}1=\sum_{m\in(\mathbb{Z}/M\mathbb{Z})^{\times}}\sum_{C\in Cl(K,\mathfrak{a})}\sum_{\mathfrak{p}\in\mathscr{P}(x,m,M,C,\mathfrak{a})}1+O(1)=\frac{x}{\log x}+o\left(\frac{x}{\log x}\right)
    $$
    as $x\to +\infty$. This is the classical Prime Ideal Theorem.
\end{remark}

\subsection{Proofs of Theorem \ref{2.1} and Theorem \ref{2.2}}
\label{section 2.2}
We introduce some notations for convenience. Let $K$ be a number field. We denote by $\mathbb{A}_K$ (resp. $\mathbb{A}_{K}^{\times}$, $\mathbb{C}_{K}=\mathbb{A}_{K}^{\times}/K^{\times}$) the adele ring (resp. the idele group, the idele class group) of $K$. For a place $v$ of $K$, we denote by $K_v$ the completion of $K$ at $v$, with $\mathfrak{o_v}$ the integer ring when $v$ is finite. A Hecke character of $K$ is a continuous character of $\mathbb{C}_{K}$, or equivalently, a continuous character of $\mathbb{A}_{K}^{\times}$ which is trivial on $K^{\times}$. 

Let $\psi$ be a unitary Hecke character of $K$, which can be decomposed as $\psi=\otimes_{v}'\psi_v$ for some continuous character $\psi_v$ of $K_{v}^{\times}$ with respect to each place $v$ of $K$. Here $\otimes_{v}'\psi_v$ is the character sending $(x_v)_v\in\mathbb{A}_{K}^{\times}$ to $\prod_v\psi_v(x_v)$, which is well defined because $\psi_v|_{\mathfrak{o}_{v}^{\times}}$ is trivial for almost all finite places $v$. For a finite place $v$, say $\psi_v$ is unramified if $\psi_v|_{\mathfrak{o}_{v}^{\times}}$ is trivial. Let $S$ be the finite set consisting of infinite places and ramified places. The Hecke $L$-function associated with $\psi$ is defined as
$$
L(\psi,s)=\prod_{v\notin S}(1-\psi_{v}(\pi_v) (Nv)^{-s})^{-1},
$$
where $\pi_v$ is a uniformizer of $K_v$ and $N$ is the norm of ideals.

\begin{proof}[Proof of Theorem \ref{2.1}]
Using the orthogonality relations for characters, we have
$$
\begin{aligned}
\sum_{\mathfrak{p}\in\mathscr{P}(x,m,M,C,\mathfrak{a})}1&=\frac{1}{\varphi(M)}\sum_{\mathfrak{p}\in\mathscr{P}(x,C,\mathfrak{a})}\sum_{\phi\in\widehat{(\mathbb{Z}/M\mathbb{Z})^{\times}}}\overline{\phi}(m)\phi(N(\mathfrak{p}))\\
&=\frac{1}{\varphi(M)  h_\mathfrak{a}}\sum_{\phi\in\widehat{(\mathbb{Z}/M\mathbb{Z})^{\times}}}\overline{\phi}(m)\sum_{N(\mathfrak{p})\leq x}\phi(N(\mathfrak{p}))\sum_{\psi\in\widehat{Cl(K,\mathfrak{a})}}\overline{\psi}(C)\psi([\mathfrak{p}])\\
&=\frac{1}{\varphi(M)  h_\mathfrak{a}}\sum_{\phi\in\widehat{(\mathbb{Z}/M\mathbb{Z})^{\times}}}\sum_{\psi\in\widehat{Cl(K,\mathfrak{a})}}\overline{\phi}(m)\overline{\psi}(C)\sum_{N(\mathfrak{p})\leq x}\phi(N(\mathfrak{p}))\psi([\mathfrak{p}]).
\end{aligned}
$$
Here, if $(N(\mathfrak{p}),M)>1$, take $\phi(N(\mathfrak{p}))=0$ and if $\mathfrak{p}$ divides $\mathfrak{a}$, take $\psi([\mathfrak{p}])=0$. 

For fixed $\phi\in\widehat{(\mathbb{Z}/M\mathbb{Z})^{\times}}, \psi\in\widehat{Cl(K,\mathfrak{a})}$, define an $L$-function
$$
L(\phi,\psi,s)=\prod_{\mathfrak{p}}(1-\phi(N(\mathfrak{p}))\psi([\mathfrak{p}])(N(\mathfrak{p}))^{-s})^{-1}.
$$
We will show $L(\phi,\psi,s)$ is a Hecke $L$-function up to finitely many Euler factors.

First we lift $\phi$ to a unitary Hecke character of $K$. Define $\widetilde{\phi}$ to be the following composite map
$$
\mathbb{C}_K\stackrel{\widetilde{N_{K/\mathbb{Q}}}}{\longrightarrow}\mathbb{C}_{\mathbb{Q}}\longrightarrow\mathbb{Q}^{\times}\backslash\mathbb{A}_{\mathbb{Q}}^{\times}/\mathbb{R}_{>0}\cong\prod_{p}\mathbb{Z}_{p}^{\times}\longrightarrow(\mathbb{Z}/M\mathbb{Z})^{\times}\stackrel{\phi^{-1}}{\longrightarrow}\mathbb{C}^{\times}.
$$
Here the isomorphism is due to the strong approximation of the idele group $\mathbb{A}_{K}^{\times}$ and $\widetilde{N_{K/\mathbb{Q}}}$ is induced by the global norm map
$$
N_{K/\mathbb{Q}}:\mathbb{A}_K\longrightarrow\mathbb{A}_{\mathbb{Q}}, (x_w)_w\longmapsto\left(\prod_{w|v}N_{w|v}(x_w)\right)_v,
$$
where $N_{w|v}: K_{w}\longrightarrow \mathbb{Q}_v$ is the norm map of the corresponding local fields. Then for any finite place $v$ of $K$ satisfying $(Nv,M)=1$, the character $\widetilde{\phi}_v$ is unramified and $\widetilde{\phi}_v(\pi_v)=\phi(Nv)$.

Next we lift $\psi$ to a unitary Hecke character $\widetilde{\psi}$ of $K$. Recall we have an isomorphism $\mathbb{C}_{K}/\overline{U(\mathfrak{a})}\cong Cl(K,\mathfrak{a})$, where $\overline{U(\mathfrak{a})}$ is the image of $U(\mathfrak{a})=\prod_{v}U_{v}(\mathfrak{a})\subset \mathbb{A}_{K}^{\times}$ in $\mathbb{C}_{K}$. Here $U_{v}(\mathfrak{a})=\ker(\mathfrak{o}_{v}^{\times}\to(\mathfrak{o}_{v}/\mathfrak{a}\mathfrak{o}_v)^{\times})$ if $v$ is a finite place, $U_{v}(\mathfrak{a})=\mathbb{C}^{\times}$ if $v$ is a complex place and $U_{v}(\mathfrak{a})=\mathbb{R}_{>0}$ if $v$ is a real place. Define $\widetilde{\psi}$ as the inflation of $\psi$. Then for any finite place $v$ not dividing $\mathfrak{a}$ of $K$, the character $\widetilde{\psi}_v$ is unramified and $\widetilde{\psi}_v(\pi_v)=\psi([v])$ .

The discussion above shows that $L(\phi,\psi,s)$ is equal to the Hecke $L$-function $L(\widetilde{\phi}\cdot\widetilde{\psi},s)$ up to finitely many Euler factors. Thus they have the same analytic properties. 

It is well known that the Hecke $L$-function $L(\widetilde{\phi}\cdot \widetilde{\psi},s)$ converges for $\Re(s)>1$ and extends to a meromorphic function on the whole plane. Moreover, it has poles if and only if $\widetilde{\phi}\cdot \widetilde{\psi}=|\cdot|_{\mathbb{A}_{K}}^{\lambda_0}$ for some pure imaginary number $\lambda_0$. In this case, take $v$ an infinite place and we have $\widetilde{\phi}_v\cdot \widetilde{\psi}_v=|\cdot|_{v}^{\lambda_0}$. However, by definition $\widetilde{\phi}_v\cdot \widetilde{\psi}_v$ is trivial when $v$ is a complex place and $\widetilde{\phi}_v\cdot \widetilde{\psi}_v|_{\mathbb{R}_{>0}}$ is trivial when $v$ is a real place, both forcing $\lambda_0=0$. Therefore, $L(\widetilde{\phi}\cdot \widetilde{\psi},s)$ has poles if and only if $\widetilde{\phi}\cdot \widetilde{\psi}$ is trivial, in which case $L(\widetilde{\phi}\cdot \widetilde{\psi},s)=\zeta_{K}(s)$. The last condition is equivalent to $\phi(N(\mathfrak{p}))\psi([\mathfrak{p}])=1$
for almost all prime ideals $\mathfrak{p}$ of $\mathcal{O}_K$ because a Hecke character is determined by its local multiplicative characters of almost all places (more strongly, of places with density strictly $>\frac{1}{2}$).

In conclusion, the $L$-function $L(\phi,\psi,s)$ converges for $\Re(s)>1$, extends to a meromorphic function on $\Re(s)\geq1$ and has neither zero nor pole except possibly a simple pole at $s=1$, and this occurs if and only if $\phi(N(\mathfrak{p}))\psi([\mathfrak{p}])=1$
for almost all prime ideals $\mathfrak{p}$ of $\mathcal{O}_K$. Apply \cite[Chapter 1, Appendix, Theorem 1]{Serre1989AbelianLR} by taking $G=(\mathbb{Z}/M\mathbb{Z})^{\times}\times Cl(K,\mathfrak{a})$ and $x_\mathfrak{p}=(N(\mathfrak{p}),[\mathfrak{p}])$ (with notations loc. cit.), we get 
$$
\sum_{\mathfrak{p}\in\mathscr{P}(x,m,M,C,\mathfrak{a})}1=\frac{A(m,M,C,\mathfrak{a})  x}{h_{\mathfrak{a}} \varphi(M) \log x }+o\left(\frac{x}{\log x}\right)
$$
as $x\to +\infty$. The result follows.
\end{proof}

\begin{proof}[Proof of Theorem \ref{2.2}]
Note that $\{x\mapsto x^n,n\in\mathbb{Z}\}$ are all characters of the circle $S^1$. By \cite[Chapter 1, Appendix, Corollary 1]{Serre1989AbelianLR} and Theorem \ref{2.1} above, it suffices to show
$$
\sum_{\mathfrak{p}\in\mathscr{P}(x,m,M,C,\mathfrak{a})}\lambda^{n}(\mathfrak{p})=o\left(\frac{x}{\log x}\right), n\neq 0
$$
as $x\to +\infty$. Using orthogonal relations of characters, we have
$$
\begin{aligned}
\sum_{\mathfrak{p}\in\mathscr{P}(x,m,M,C,\mathfrak{a})}\lambda^{n}(\mathfrak{p})&=\frac{1}{\varphi(M)}\sum_{\mathfrak{p}\in\mathscr{P}(x,C,\mathfrak{a})}\sum_{\phi\in\widehat{(\mathbb{Z}/M\mathbb{Z})^{\times}}}\overline{\phi}(m)\phi(N(\mathfrak{p}))\lambda^{n}(\mathfrak{p})\\
&=\frac{1}{\varphi(M)  h_\mathfrak{a}}\sum_{\phi\in\widehat{(\mathbb{Z}/M\mathbb{Z})^{\times}}}\overline{\phi}(m)\sum_{N(\mathfrak{p})\leq x}\phi(N(\mathfrak{p}))\sum_{\psi\in\widehat{Cl(K,\mathfrak{a})}}\overline{\psi}(C)\psi([\mathfrak{p}])\lambda^{n}(\mathfrak{p})\\
&=\frac{1}{\varphi(M)  h_\mathfrak{a}}\sum_{\phi\in\widehat{(\mathbb{Z}/M\mathbb{Z})^{\times}}}\sum_{\psi\in\widehat{Cl(K,\mathfrak{a})}}\overline{\phi}(m)\overline{\psi}(C)\sum_{N(\mathfrak{p})\leq x}\phi(N(\mathfrak{p}))\psi([\mathfrak{p}])\lambda^{n}(\mathfrak{p}).
\end{aligned}
$$
For $\phi\in\widehat{(\mathbb{Z}/M\mathbb{Z})^{\times}}$, $\psi\in\widehat{Cl(K,\mathfrak{a})}$ and $n\neq 0$, define an $L$-function
$$
L(\phi,\psi,\lambda^n,s)=\prod_{\mathfrak{p}}(1-\phi(N(\mathfrak{p}))\psi([\mathfrak{p}])\lambda^{n}(\mathfrak{p})(N(\mathfrak{p}))^{-s})^{-1}.
$$
We will show that up to finitely many Euler factors $L(\phi,\psi,\lambda^n,s)$ is a Hecke $L$-function with no pole on $\Re(s)=1$, and the result follows by \cite[Chapter 1, Appendix, Theorem 1]{Serre1989AbelianLR}, taking $G=(\mathbb{Z}/M\mathbb{Z})^{\times}\times Cl(K,\mathfrak{a})\times S^1$ and $x_\mathfrak{p}=(N(\mathfrak{p}),[\mathfrak{p}],\lambda(\mathfrak{p}))$.

Lift $\phi$, $\psi$ to unitary Hecke characters $\widetilde{\phi}$, $\widetilde{\psi}$ of $K$ as above. We lift $\lambda$ to a unitary Hecke character $\widetilde{\lambda}=\otimes_{v}^{'} \widetilde{\lambda}_v$ of $K$: for $v$ finite place, we take $\widetilde{\lambda}_v$ to be unramified and define $\widetilde{\lambda}_v(\pi_v)=\lambda(v)$; for $v=\infty$ complex place, define $\widetilde{\lambda}_{\infty}(z)=\left(\frac{z}{|z|}\right)^{-w}$ for $z\in \mathbb{C}^{\times}$. Then $L(\phi,\psi,\lambda^{n},s)$ is equal to the Hecke $L$-function $L(\widetilde{\phi}\cdot\widetilde{\psi} \cdot\widetilde{\lambda}^n,s)$ up to finitely many Euler factors. It has poles if and only if $\widetilde{\phi}\cdot \widetilde{\psi}\cdot{\widetilde{\lambda}}^{n}=|\cdot|_{\mathbb{A}_{K}}^{\lambda_0}$ for some pure imaginary number $\lambda_0$, in which case we have $|x|_{\infty}^{\lambda_0}=\widetilde{\lambda}_{\infty}^n(x)=(\frac{x}{|x|})^{-wn}$ for all $x\in \mathbb{C}^{\times}$. This is impossible.
\end{proof}

\begin{remark}
\label{remark2.4}
    Sometimes it is useful to consider the argument-equidistribution of primary generators instead of prime ideals. For example, we take $K=\mathbb{Q}(\sqrt{-1})$. Instead of the argument homomorphism $\lambda$ we used above, we define a new homomorphism $\Lambda:P(K)\rightarrow S^{1}$ as follows: for a nonzero prime ideal $\mathfrak{p}$ prime to $2$, there exists a unique generator $\alpha$ of $\mathfrak{p}$ satisfying $\alpha \equiv 1 \mod(1+\sqrt{-1})^3$, and we define $\Lambda(\mathfrak{p})=\frac{\alpha}{|\alpha|}$; we define $\Lambda\left((1+\sqrt{-1})\right)=\frac{1+\sqrt{-1}}{\sqrt{2}}$. To obtain the argument-equidistribution in terms of $\Lambda$ in the sense of Theorem \ref{2.2}, we can argue as in the proof of Theorem \ref{2.2} above. The only difference is how to lift $\Lambda$ to a Hecke character $\widetilde{\Lambda}=\otimes_{v}^{'} \widetilde{\Lambda}_v$ of $K$: for $v=\infty$ complex place, define $\widetilde{\Lambda}_{\infty}(z)=\left(\frac{z}{|z|}\right)^{-1}$ for $z\in \mathbb{C}^{\times}$; for $v$ finite place prime to $2$, we take $\widetilde{\Lambda}_v$ to be unramified and define $\widetilde{\Lambda}_v(\pi_v)=\Lambda(v)$; for $v=(1+\sqrt{-1})$, we let $\widetilde{\Lambda}_v|_{\mathfrak{o}_{v}^{\times}}$ trivial on $1+\pi_v^3\mathfrak{o}_{v}$ and define $\widetilde{\Lambda}_v(\sqrt{-1})=\sqrt{-1},~\widetilde{\Lambda}_v(1+\sqrt{-1})=\Lambda(v)$. This determines $\widetilde{\Lambda}_v$ because the image of $\sqrt{-1}$ is a generator of $\mathfrak{o}_{v}^{\times}/(1+\pi_v^3\mathfrak{o}_{v})\cong\mathbb{Z}/4\mathbb{Z}$. One checks $\widetilde{\Lambda}=\otimes_{v}^{'} \widetilde{\Lambda}_v$ constructed above is indeed a Hecke character of $K$. Therefore, we obtain the argument-equidistribution analogous to Theorem \ref{2.2}, in terms of $\Lambda$. This motivates us to consider the counting functions in Section \ref{4.2}.
\end{remark}

\subsection{View of $\check{\mathrm{C}}$ebotarev density theory}
\label{section 2.3}
In this section, we reinterpret Theorem \ref{2.1} from the view of $\check{\mathrm{C}}$ebotarev density theory. 

Take $K,m,M,\mathfrak{a},C$ as in Theorem \ref{2.1}. We denote by $K^{ab}$ the maximal abelian extension of $K$. By the global class field theory, we have an Artin map
$$
r_K:\mathbb{C}_K\longrightarrow {\rm{Gal}}(K^{ab}/K) 
$$
which induces an isomorphism
$$
\widetilde{r_K}:Cl(K,\mathfrak{a})\longrightarrow {\rm{Gal}}(K(\mathfrak{a})/K),
$$
where $K(\mathfrak{a})$ is the fixed field of $r_K(\overline{U(\mathfrak{a})})$ under the Galois correspondence (see the proof of Theorem \ref{2.1} for the definition of $U(\mathfrak{a})$). For a prime ideal $\mathfrak{p}$ not dividing $\mathfrak{a}$, the isomorphism takes $\mathfrak{p}$ to $\left(\frac{K(\mathfrak{a})/K}{\mathfrak{p}}\right)$, the Frobenius element of $\mathfrak{p}$. Thus the condition $[\mathfrak{p}]=C$ can be translated into $\left(\frac{K(\mathfrak{a})/K}{\mathfrak{p}}\right)=\widetilde{r_K}(C)$. 

Let $\zeta_M$ be a primitive $M$-th root of unity, $\mathbb{Q}(\zeta_M)$ the $M$-th cyclotomic field and 
$$
\widetilde{r_{\mathbb{Q}}}:(\mathbb{Z}/M\mathbb{Z})^{\times}\longrightarrow {\rm{Gal}}(\mathbb{Q}(\zeta_M)/\mathbb{Q})
$$
the canonical isomorphism. Let $L_1=K\cap\mathbb{Q}(\zeta_M)$. We have a canonical isomorphism
$$
{\rm{Gal}}(K(\zeta_M)/K)\stackrel{\rm{Res}}{\longrightarrow}{\rm{Gal}}(\mathbb{Q}(\zeta_M)/L_1).
$$
Take a prime ideal $\mathfrak{p}$ of $K$ such that $(N(\mathfrak{p}),M)=1$. Let $\mathfrak{p}\cap \mathbb{Q}=p\mathbb{Z}$. By the properties of Frobenius under base change, we have:
$$
\left(\frac{K(\zeta_M)/K}{\mathfrak{p}}\right)\bigg|_{\mathbb{Q}(\zeta_M)}=\left(\frac{\mathbb{Q}(\zeta_M)/\mathbb{Q}}{p}\right)^{f(\mathfrak{p}|p)}=\widetilde{r_{\mathbb{Q}}}(N(\mathfrak{p})).
$$
If $\widetilde{r_{\mathbb{Q}}}(m)\notin {\rm{Gal}}(\mathbb{Q}(\zeta_M)/L_1)$, there is no $\mathfrak{p}$ satisfying $N(\mathfrak{p})\equiv m\pmod{M}$. Suppose $\widetilde{r_{\mathbb{Q}}}(m)\in {\rm{Gal}}(\mathbb{Q}(\zeta_M)/L_1)$ and view it as an element of ${\rm{Gal}}(K(\zeta_M)/K)$. Then the condition $N(\mathfrak{p})\equiv m\pmod{M}$ can be translated into $\left(\frac{K(\zeta_M)/K}{\mathfrak{p}}\right)=\widetilde{r_{\mathbb{Q}}}(m)$.

Take $L_2=K(\mathfrak{a})\cap K(\zeta_M)$, $\widetilde{L_2}=L_2\cap \mathbb{Q}(\zeta_M)=K(\mathfrak{a})\cap \mathbb{Q}(\zeta_M)$ and $L_3=K(\mathfrak{a})  K(\zeta_M)$.
\[\begin{tikzcd}
	&& {L_3} \\
	{K(\mathfrak{a})} && {K(\zeta_M)} \\
	{L_2} &&& {\mathbb{Q}(\zeta_M)} \\
	K & {\widetilde{L_2}} \\
	& {L_1} \\
	& {\mathbb{Q}}
	\arrow[no head, from=1-3, to=2-3]
	\arrow[no head, from=2-1, to=1-3]
	\arrow[no head, from=2-1, to=3-1]
	\arrow[no head, from=2-3, to=3-1]
	\arrow[no head, from=2-3, to=3-4]
	\arrow[no head, from=3-1, to=4-1]
	\arrow[no head, from=3-1, to=4-2]
	\arrow[no head, from=4-1, to=2-3]
	\arrow[no head, from=4-1, to=5-2]
	\arrow[no head, from=4-2, to=3-4]
	\arrow[no head, from=4-2, to=5-2]
	\arrow[no head, from=5-2, to=3-4]
	\arrow[no head, from=5-2, to=6-2]
\end{tikzcd}\]
Note that $\left(\frac{K(\mathfrak{a})/K}{\mathfrak{p}}\right)\bigg|_{L_2}=\left(\frac{L_2/K}{\mathfrak{p}}\right)=\left(\frac{K(\zeta_M)/K}{\mathfrak{p}}\right)\bigg|_{L_2}$. Thus if $\widetilde{r_K}(C)|_{L_2}\neq \widetilde{r_{\mathbb{Q}}}(m)|_{L_2}$, there does not exist $\mathfrak{p}$ satisfying $[\mathfrak{p}]=C$ and $N(\mathfrak{p})\equiv m\pmod{M}$ simultaneously. Suppose $\widetilde{r_K}(C)|_{L_2}= \widetilde{r_{\mathbb{Q}}}(m)|_{L_2}$. Then $\widetilde{r_K}(C)$ and $\widetilde{r_{\mathbb{Q}}}(m)$ uniquely determine an element $\sigma_{m,C}\in {\rm{Gal}}(L_3/K)$. Hence the condition $[\mathfrak{p}]=C$, $N(\mathfrak{p})\equiv m\pmod{M}$ can be translated into $\sigma_{m,C}=\left(\frac{L_3/K}{\mathfrak{p}}\right)$. By $\check{\mathrm{C}}$ebotarev density theory, we have
$$
\begin{aligned}
\sum_{\mathfrak{p}\in\mathscr{P}(x,m,M,C,\mathfrak{a})}1&=\frac{1}{[L_3:K]} \frac{x}{\log x}+o\left(\frac{x}{\log x}\right)\\
&=\frac{[L_2:K]}{[K(\zeta_M):K]  h(\mathfrak{a})} \frac{x}{\log x}+o\left(\frac{x}{\log x}\right)\\
&=\frac{[\widetilde{L_2}:\mathbb{Q}]}{\varphi(M)  h(\mathfrak{a})} \frac{x}{\log x}+o\left(\frac{x}{\log x}\right).
\end{aligned}
$$
as $x\to +\infty$. 

Now we explain that this is consistent with Theorem \ref{2.1}. Recall
$$
A(m,M,C,\mathfrak{a})=\sum_{\phi,\psi}\phi(m)\psi(C),
$$
the sum ranging over $\phi\in\widehat{(\mathbb{Z}/M\mathbb{Z})^{\times}}, \psi\in\widehat{Cl(K,\mathfrak{a})}$ satisfying $\phi(N(\mathfrak{p}))\psi([\mathfrak{p}])=1$ for almost all prime ideals $\mathfrak{p}$ of $\mathcal{O}_K$. Identify $\phi$ (resp. $\psi$) with character of ${\rm{Gal}}(\mathbb{Q}(\zeta_M)/\mathbb{Q})$ (resp. ${\rm{Gal}}(K(\mathfrak{a})/K)$) through $\widetilde{r_\mathbb{Q}}$ (resp. $\widetilde{r_K}$). Then 
$$
A(m,M,C,\mathfrak{a})=\sum_{\phi,\psi}\phi(\widetilde{r_\mathbb{Q}}(m))\psi(\widetilde{r_K}(C)),
$$
the sum ranging over $\phi, \psi$ satisfying $\phi\left(\left(\frac{K(\zeta_M)/K}{\mathfrak{p}}\right)\bigg|_{\mathbb{Q}(\zeta_M)}\right)\psi\left(\left(\frac{K(\mathfrak{a})/K}{\mathfrak{p}}\right)\right)=1$ for almost all prime ideals $\mathfrak{p}$ of $\mathcal{O}_K$. Define $\widetilde{\phi}$ to be the composition
$$
{\rm{Gal}}(K(\zeta_M)/K)\stackrel{\rm{Res}}{\longrightarrow}{\rm{Gal}}(\mathbb{Q}(\zeta_M)/L_1)\stackrel{\phi}{\longrightarrow}\mathbb{C}^{\times}.
$$
Then the condition on $\phi, \psi$ is equivalent to $\widetilde{\phi}|_{{\rm{Gal}}(K(\zeta_M)/L_2)},\psi|_{{\rm{Gal}}(K(\mathfrak{a})/L_2)}$ and $\widetilde{\phi} \cdot\psi$ (viewed as character of ${\rm{Gal}}(L_2/K)$) are trivial. Thus it suffices to sum over $\phi\in\widehat{{\rm{Gal}}(\widetilde{L_2}/\mathbb{Q})}$ and $A(m,M,C,\mathfrak{a})=\sum_{\phi}\phi(\widetilde{r_\mathbb{Q}}(m)|_{\widetilde{L_2}})\phi^{-1}(\widetilde{r_K}(C)|_{\widetilde{L_2}})$, which is $[\widetilde{L_2}:\mathbb{Q}]$ if $\widetilde{r_K}(C)|_{\widetilde{L_2}}= \widetilde{r_{\mathbb{Q}}}(m)|_{\widetilde{L_2}}$ and otherwise $0$.

We point out that $\check{\mathrm{C}}$ebotarev density theory is a corollary of \cite[Chapter 1, Appendix, Theorem 1]{Serre1989AbelianLR}  by taking $G={\rm{Gal}}(L/K)$ for a Galois extension $L$ over $K$ and $x_\mathfrak{p}=\left(\frac{L/K}{\mathfrak{p}}\right)$ (with notations loc. cit.) for primes $\mathfrak{p}$ unramified in $L$. The corresponding $L$-function is Artin $L$-function of $K$ up to finitely many Euler factors. One can establish a similar view for Theorem \ref{2.2} by taking $G=\rm{Gal}(L/K)\times S^1$ for a Galois extension $L$ of $K$ and $x_\mathfrak{p}=\left(\left(\frac{L/K}{\mathfrak{p}}\right), \lambda(\mathfrak{p})\right)$ for primes $\mathfrak{p}$ unramified in $L$. 

\subsection{Proof of Theorem \ref{main theorem}}
\label{section 2.4}
\begin{proof}[Proof of Theorem \ref{main theorem}]
We use the notation of Section \ref{section 2.3} with $\mathfrak{a}=\mathfrak{f}$. Recall the notation $g(\mathfrak{f}),~C\in Cl(K,\mathfrak{f}),~\mathfrak{a}_0 \in C^{-1}$ and $\lambda$ defined above Theorem \ref{Coleman}. By discussion in Section \ref{section 2.3}, we can assume $\widetilde{r_K}(C)|_{L_2}=\widetilde{r_{\mathbb{Q}}}(m)|_{L_2}$. Otherwise, there is no prime $\mathfrak{p}$ satisfying the condition. Let $\sigma_{m,C}\in{\rm{Gal}}(L_3/K)$ be the element uniquely determined by $\widetilde{r_K}(C)$ and $\widetilde{r_{\mathbb{Q}}}(m)$. The condition $[\mathfrak{p}]=C$, $N(\mathfrak{p})\equiv m\pmod{M}$ is equivalent to $\sigma_{m,C}=\left(\frac{L_3/K}{\mathfrak{p}}\right)$. Since $L_3/K$ is abelian, we may choose a fixed ideal $\mathfrak{b}$ of $\mathcal{O}_K$, divisible by $\mathfrak{f}$, such that $L_3\subset K(\mathfrak{b})$. Here $K(\mathfrak{b})$ is the fixed field of $r_K(\overline{U(\mathfrak{b})})$ under the Galois correspondence. 

Assume $[K(\mathfrak{b}):L_3]=s$ and denote the elements of ${\rm{Gal}}(K(\mathfrak{b})/K)$ extending $\sigma_{m,C}$ by $\sigma_i,1\leq i\leq s$. Put $C_i=\overline{r_K}^{-1}(\sigma_i)$, where $\overline{r_K}$ is the isomorphism
$$
\overline{r_K}:Cl(K,\mathfrak{b})\xrightarrow{\sim} {\rm{Gal}}(K(\mathfrak{b})/K).
$$
Then the condition $\sigma_{m,C}=\left(\frac{L_3/K}{\mathfrak{p}}\right)$ is equivalent to $[\mathfrak{p}]\in\{C_i\mid 1\leq i\leq s\}$. For each $i$, choose an ideal
$\mathfrak{a}_{0,i}\in C_i^{-1}$. For $\mathfrak{p}\in C_i$, we can find some $\xi_{\mathfrak{p},i}\in K$ satisfying $(\xi_{\mathfrak{p},i})=\mathfrak{p}\mathfrak{a}_{0,i}$ and $\xi_{\mathfrak{p},i}\equiv1\pmod{\mathfrak{b}}$ and we write $\lambda_i(\mathfrak{p})=\left(\frac{\xi_{\mathfrak{p},i}}{|\xi_{\mathfrak{p},i}|}\right)^{g(\mathfrak{b})}$. Put $r=\frac{g(\mathfrak{f})}{g(\mathfrak b)}\in\mathbb{Z}$. Since $\mathfrak{a}_0$ and $\mathfrak{a}_{0,i}$ represent the same class modulo $\mathfrak{f}$, there exists some $\delta_i\equiv1\pmod{\mathfrak{f}}$ such that $\mathfrak{a}_0=(\delta_i)\mathfrak{a}_{0,i}$. Then $\lambda(\mathfrak{p})=\left(\frac{\delta_i}{|\delta_i|}\right)^{g(\mathfrak{f})}\lambda_i(\mathfrak{p})^r$. Write $\left(\frac{\delta_i}{|\delta_i|}\right)^{g(\mathfrak{f})}=e^{it_i}$. Then the angular condition $arg(\lambda(\mathfrak{p}))\in[\varphi_1,\varphi_2]$ is equivalent to $arg(\lambda_i(\mathfrak{p}))\in\bigcup_{j=1}^rJ_{i,j}$, where
\[
J_{i,j}=\left[\frac{\varphi_1-t_i+2\pi j}{r},\frac{\varphi_2-t_i+2\pi j}{r}\right]\pmod{2\pi}.
\]
Applying Theorem \ref{Coleman} to modulus
$\mathfrak b$, class $C_i$ and arcs $J_{i,j}$ for $1\le i\le s,~1\le j\le r$, we obtain
\begin{align*}
    \sum_{\mathfrak{p}\in P(x,y,C,\varphi_1,\varphi_2,m,M)}1&=\frac{s(\varphi_2-\varphi_1)  y}{2\pi  h_{\mathfrak{b}} \log x }\left(1+O\left(\frac{1}{\log x}\right)\right)\\
    &=\frac{(\varphi_2-\varphi_1)  y}{2\pi [L_3:K]\log x }\left(1+O\left(\frac{1}{\log x}\right)\right)\\
    &=\frac{A(m,M,C,\mathfrak{f}) (\varphi_2-\varphi_1)  y}{2\pi  h_{\mathfrak{f}} \varphi(M) \log x }\left(1+O\left(\frac{1}{\log x}\right)\right).
\end{align*}
\end{proof}

\section{An application}
\label{pf of application}
As an application of Theorem \ref{main theorem}, we generalize C. Elsholtz and G. Harman's work \cite{MR3467390} on the conjectures of T. Ordowski and Z.-W. Sun \cite{sun2017conjectures}. 

For better illustration, we briefly recall some fundamental concepts regarding binary quadratic forms. A binary quadratic form $Q(x,y)=ax^2+bxy+cy^2$ with integer coefficients is called \textit{primitive} if $\gcd(a,b,c)=1$. We say that a quadratic form can represent a prime number $p$ if there exist integers $x_p, y_p$ such that $p=Q(x_p,y_p)$. A quadratic form is \textit{positive definite} if $Q(x,y)>0$ for all $(x,y) \neq (0,0)$, which is equivalent to $a>0$ and its discriminant $\Delta = b^2-4ac < 0$. A positive definite quadratic form can represent infinitely many prime numbers only if it is primitive. 

Elsholtz and Harman \cite{MR3467390} considered the limits of coordinate-wise defined functions over all representable primes. Our Theorem \ref{main} generalizes their work to the set of representable primes satisfying additional congruence conditions. We follow the proof of Theorem 1.1 in \cite{MR3467390}.

\begin{proof}[Proof of Theorem \ref{main}]
Let $\Delta=-D$ and $K=\mathbb{Q}(\sqrt{\Delta})$. We first assume that $\Delta$ is a fundamental discriminant and $\Delta<-4$. Let
$$\mathfrak{c} = \left[a, \frac{1}{2}(b+\sqrt{-D})\right], \quad \mathfrak{d} = \left[a, \frac{1}{2}(b-\sqrt{-D})\right]$$
be two ideals in $\mathbb{Q}(\sqrt{\Delta})$. If $\mathfrak{c}$ belongs to an ideal class $C \in Cl(K)$, then $\mathfrak{d}$ belongs to $C^{-1}$. According to \cite[Theorem 2.2]{Coleman1990TheDO}, which is based on \cite[Theorem 2.1]{Knap}, there exists a 1-1 correspondence between the prime ideals $\mathfrak{p}\in C$ with norm $N(\mathfrak{p})=p$ and pairs of coprime integers $(x_p,y_p)$ satisfying $Q(x_p,y_p)=p$ with $x_p>0$ or $x_p=0, y_p=1$. Moreover, for two correspondents $\mathfrak{p}$ and $(x_p,y_p)$, we have $\mathfrak{pd}=(\xi_\mathfrak{p})$ where $\xi_\mathfrak{p}=x_p a + y_p \left( \frac{b - \sqrt{-D}}{2} \right)$. Taking the absolute value, we obtain $|\xi_{\mathfrak{p}}|^2 = a p$. 

Writing $r = \sqrt{p}$ and defining the argument $\theta = -\arg(\xi_{\mathfrak{p}})$, we can invert this relation to express the coordinates $(x_p, y_p)$ in terms of $r$ and $\theta$
$$
x_p = \frac{r}{\sqrt{aD}} (\sqrt{D} \cos\theta - b \sin\theta), \quad y_p = 2a \frac{r}{\sqrt{aD}} \sin\theta.
$$
The condition $x_p > y_p>0$ imposed on the coordinates restricts the argument $\theta$ to a specific sector $(0, \beta)$ in the complex plane, where $\beta$ is determined by the coefficients of $Q(x,y)$ as defined in Theorem \ref{main}. Under the 1-1 correspondence above, this requires the corresponding prime ideal $\mathfrak{p}$ satisfying $arg(\lambda(\mathfrak{p}))\in (-g\beta,0)$ with $g=2$ and $\mathfrak{a}_0=\mathfrak{d}$ (recall the notation $g$, $\mathfrak{a}_0$ and $\lambda$ defined above Theorem \ref{Coleman}).

Under the established correspondence, requiring the prime $p \equiv m \pmod M$ is strictly equivalent to imposing the norm condition $N(\mathfrak{p}) \equiv m \pmod M$ on the corresponding prime ideals. This exact algebraic translation allows us to apply Theorem \ref{main theorem} directly, setting $\mathfrak{f}=\mathcal{O}_K$ and $C\in Cl(K)$ as above (note that the assumption that $P_{m,M}\cap\{Q(x,y)\mid x,y \in\mathbb{Z}\}$ has positive density ensures $A(m,M,C,\mathcal{O}_K)\neq 0$). We dissect the region $0 < r \le \sqrt{N}$ and $0 \le \theta \le \beta$ into fine polar boxes as in the proof of Theorem 1.1 in \cite{MR3467390}. According to our Theorem \ref{main theorem}, the prime ideals in the ideal class $C$ satisfying $N(\mathfrak{p}) \equiv m \pmod M$ are asymptotically equidistributed with respect to their arguments within this sector. Converting the summations over the polar boxes into Riemann integrals as $N \to \infty$, we have
\begin{equation*}
\lim_{N\rightarrow\infty}\frac{\sum\limits_{p\leq N,p\in P_{m,M}}~\sum\limits_{(x_p,y_p)\in \Omega_p}x_p^{k}}{\sum\limits_{p\leq N,p\in P_{m,M}}~\sum\limits_{(x_p,y_p)\in \Omega_p}y_p^{k}} = \frac{\int_{0}^{\beta} (\sqrt{D} \cos \theta - b \sin \theta)^k \dif\theta}{\int_{0}^{\beta} (2a \sin \theta)^k \dif\theta}=\frac{\int_0^\beta s(k, \theta) \dif\theta}{\int_0^\beta t(k, \theta) \dif\theta}.
\end{equation*}

Now assume that $\Delta=-4$ or $-3$. Then $K=\mathbb{Q}(\sqrt{-1})$ or $\mathbb{Q}(\sqrt{-3})$, which has $4$ or $6$ units respectively. In this case, the correspondence from integer pairs $(x_p,y_p)$ to prime ideals $\mathfrak{p}$ may not be one-to-one, even restricting to $x_p > y_p>0$: for prime ideal $\mathfrak{p}\in C$, the generator $\xi_{\mathfrak{p}}$ of the principal ideal $\mathfrak{pd}$ is only unique up to multiplication by units, and the sector $(0,\beta)$ can contain more than one arguments $\theta = -\arg(\xi_{\mathfrak{p}})$ if $\beta>\frac{\pi}{2}$ and $\Delta=-4$, or $\beta>\frac{\pi}{3}$ and $\Delta=-3$ . For example, consider $Q(x,y)=x^2-4xy+5y^2$ with $\Delta=-4$ and $\beta=\arctan(-1)=\frac{3\pi}{4}$. The prime ideal $\mathfrak{p}=(2-i)=(-1-2i)$ corresponds to two integer pairs $(4,1)$ and $(3,2)$, with corresponding arguments $\theta=\arctan(\frac{1}{2})$ and $\theta=\arctan(-2)$ respectively. To resolve this problem, we can partition the interval $(0,\beta)$ into disjoint subintervals with length less than $\frac{\pi}{2}$ if $\Delta=-4$ and $\frac{\pi}{3}$ if $\Delta=-3$. For each subinterval $J$, there is a 1-1 correspondence between pairs of coprime integers $(x_p,y_p)$ with $Q(x_p,y_p)=p$ satisfying $\theta\in J$, and the prime ideals $\mathfrak{p}$ with norm $N(\mathfrak{p})=p$ satisfying $arg(\lambda(\mathfrak{p}))\in -gJ\pmod{2\pi}$, where $g=4$ if $\Delta=-4$ and $g=6$ if $\Delta=-3$. Then we can apply Theorem \ref{main theorem} to each subinterval $J$ and obtain the same result.

Finally, suppose that $\Delta=f^2\Delta_0$, where $\Delta_0$ is a fundamental discriminant and $f>1$. Take $\mathfrak d_f=\left[a,\frac{b-f\sqrt{\Delta_0}}2\right]$ and $\mathcal{O}_f=\mathbb{Z}+f\mathcal{O}_K$. In this case, the correspondence is between coprime integers $(x_p,y_p)$ and prime ideals $\mathfrak{p}$ of the order $\mathcal{O}_f$ satisfying $\mathfrak{p}\in [\mathfrak d_f]^{-1}\in \mathrm{Pic}(\mathcal{O}_f)$ and $N(\mathfrak{p})=p$. Here we exclude a finite number of primes $p$ dividing $f$. After choosing a suitable ideal $\mathfrak{a}$ (depending on $f$) of $\mathcal{O}_K$, prime ideals of $\mathcal{O}_f$ in $[\mathfrak d_f]^{-1}$ with norm $p$ correspond to prime ideals of $\mathcal{O}_K$ with norm $p$ in certain ideal classes $C_1,C_2,...,C_r\in Cl(K,\mathfrak{a})$. We refer readers to the end of the proof of Theorem 1.1 in \cite{MR3467390} for an illustrative example. Then apply Theorem 1.2 to classes $C_1,C_2,...,C_r$ respectively and then sum the resulting asymptotic formulas. All the arithmetic coefficients in the numerator and denominator cancel out and the same limiting ratio holds.
\end{proof}

\begin{example}
Consider the binary quadratic form $Q(x,y)=x^2+y^2$, whose discriminant is $\Delta = -4$. A prime $p$ can be represented as $p=x_p^2+y_p^2$ with $x_p > y_p > 0$ if and only if $p \equiv 1 \pmod 4$, in which case, such representation is unique. Applying Theorem \ref{main}, the asymptotic coordinate ratios over all representable primes in the congruence class $p \equiv 1 \pmod 8$ and $p \equiv 5 \pmod 8$ are
\small
$$\lim_{N\to\infty}\frac{\sum\limits_{p\leq N,p\equiv1\pmod{8}}x_p}{\sum\limits_{p\leq N,p\equiv1\pmod{8}}y_p}=\lim_{N\to\infty}\frac{\sum\limits_{p\leq N,p\equiv5\pmod{8}}x_p}{\sum\limits_{p\leq N,p\equiv5\pmod{8}}y_p}=\lim_{N\to\infty}\frac{\sum\limits_{p\leq N,p\equiv1\pmod{4}}x_p}{\sum\limits_{p\leq N,p\equiv1\pmod{4}}y_p}=1+\sqrt{2}.$$
This invariance naturally motivates us to investigate the secondary-term discrepancies between different residue classes, which we explore in Section 4.1.
\end{example}

We can further generalize this result from monomials to arbitrary bivariate polynomial functions of the coordinates. For a bivariate polynomial $f(x,y)$ of degree $n$, let $\tilde{f}(x,y)$ denote its homogeneous component of degree $n$. The following theorem demonstrates that the highest degree homogeneous part determines the asymptotic behavior.

\begin{theorem}
\label{polynomial}
Use the same notation as in Theorem 1.3. Suppose that $f(x,y)$ and $g(x,y)$ are two bivariate polynomials over $\mathbb{R}$ with $\deg f(x,y) = \deg g(x,y) = n\geq1$. If $P_{m,M}\cap \{Q(x,y) \mid x, y \in \mathbb{Z}\}$ has positive density in $P_{m,M}$, and $\int_0^\beta
\widetilde g(u(\theta),v(\theta))\,d\theta\neq0$, then
   $$\lim_{N\to\infty}\frac{\sum\limits_{p\leq N,p\in P_{m,M}}~\sum\limits_{(x_p,y_p)\in \Omega_p}f(x_p,y_p)}{\sum\limits_{p\leq N,p\in P_{m,M}}~\sum\limits_{(x_p,y_p)\in \Omega_p}g(x_p,y_p)}=\frac{\int_0^{\beta}{\tilde{f}(u(\theta),v(\theta))\dif\theta}}{\int_0^{\beta}{\tilde{g}(u(\theta),v(\theta))\dif\theta}},$$
   where $u(\theta)=\sqrt D\cos\theta-b\sin\theta,~v(\theta)=2a\sin\theta$.
\end{theorem}
\begin{proof}
We can decompose the polynomial $g(x,y)$ as $g(x,y) = \tilde{g}(x,y) + r_g(x,y)$, where $\tilde{g}(x,y)$ is the homogeneous polynomial of degree $n$, and $r_g(x,y)$ is a polynomial with degree at most $n-1$.

By the explicit coordinate transformation established in the proof of Theorem \ref{main}, for any represented prime $p \le N$ with $p \in P_{m,M}$, the integer coordinates $(x_p, y_p)\in\Omega_p$ are parameterized by $r = \sqrt{p}$ and the argument $\theta \in (0, \beta)$ as
$$
x_p = \frac{\sqrt{p}}{\sqrt{aD}} (\sqrt{D} \cos\theta - b \sin\theta), \quad y_p = 2a \frac{\sqrt{p}}{\sqrt{aD}} \sin\theta,
$$
where $a, b$ are fixed coefficients of the primitive positive definite quadratic form $Q(x,y)$. Thus, we have $x_p \ll p^{1/2} \le N^{1/2}$ and $y_p \ll p^{1/2} \le N^{1/2}$, where $\ll$ means that the ratio of the left-hand side to the right-hand side is bounded by a constant depending only on $a, b$ and $\Delta$. It follows that the lower degree part $r_g(x,y)$ is bounded by $O(N^{(n-1)/2})$. Hence, we have
$$
\sum\limits_{p\leq N,p\in P_{m,M}}~\sum\limits_{(x_p,y_p)\in \Omega_p} r_g(x_p, y_p) \ll \pi(N) N^{(n-1)/2} \ll \frac{N^{(n+1)/2}}{\log N},
$$
where $\pi(N)$ denotes the number of primes $p\leq N$.

On the other hand, the highest degree homogeneous part $\tilde{g}(x_p, y_p)$ can be asymptotically evaluated using partial summation and the equidistribution result from Theorem \ref{main theorem}. By dissecting the angular sector into polar boxes and integrating as $N \to \infty$, we have
$$\sum\limits_{p\leq N,p\in P_{m,M}}~\sum\limits_{(x_p,y_p)\in \Omega_p} \tilde{g}(x_p, y_p) \sim C_g \frac{N^{n/2 + 1}}{\log N},$$
where $C_g$ is a constant proportional to $A(m,M,C,\mathfrak{f})$ and the angular integral of $\tilde{g}$ over the sector $(0, \beta)$. By assumption, the constant $C_g$ is nonzero. Similar analysis holds for the polynomial  $f(x,y)$. It follows that
\begin{align*}
\lim_{N\to\infty} \frac{\sum\limits_{p\leq N,p\in P_{m,M}}~\sum\limits_{(x_p,y_p)\in \Omega_p} f(x_p, y_p)}{\sum\limits_{p\leq N,p\in P_{m,M}}~\sum\limits_{(x_p,y_p)\in \Omega_p} g(x_p, y_p)} 
&= \lim_{N\to\infty} \frac{\sum\limits_{p\leq N,p\in P_{m,M}}~\sum\limits_{(x_p,y_p)\in \Omega_p} \tilde{f}(x_p, y_p)}{\sum\limits_{p\leq N,p\in P_{m,M}}~\sum\limits_{(x_p,y_p)\in \Omega_p} \tilde{g}(x_p, y_p)}\\
&=\frac{\int_0^{\beta}{\tilde{f}(u(\theta),v(\theta))\dif\theta}}{\int_0^{\beta}{\tilde{g}(u(\theta),v(\theta))\dif\theta}}.
\end{align*}
This completes the proof.
\end{proof}

\section{Some conjectures}
A central analytical feature of our main results (Theorems 1.2, 1.3, 2.1, and 2.2) is that the leading asymptotic main terms are independent of the choice of $m$ among the reduced residue classes modulo $M$ satisfying the compatibility condition of Section \ref{section 2.3}. When the primary terms of two counting distribution functions coincide, it becomes a natural and interesting problem to investigate the dependence of the secondary terms on the congruence classes. 

Inspired by Chebyshev’s bias in  \cite{em/1048515870} and
``murmuration” in \cite{lee2024murmurations}, we investigate the difference of distributions between primes congruent to  $1 \pmod{8}$ and primes congruent to $5\pmod{8}$ and their representations by quadratic forms. We observe fascinating phenomena that invite further theoretical investigation. We list two conjectures.   

\subsection{Chebyshev’s bias}
\label{fr}
Recall $P_{m,M}=\{p\mid p\equiv m\pmod{M},p\text{ is prime}\}$ with $m,M$ positive integers and $\mathrm{gcd}(m,M)=1$. Using the same notation as in Theorem \ref{main}, for a given quadratic form $Q(x,y)$, take
$$\fr(N;M,m)=\frac{\sum\limits_{p\leq Pr(N),p\in P_{m,M}}~\sum\limits_{(x_p,y_p)\in \Omega_p}x_p}{\sum\limits_{p\leq Pr(N),p\in P_{m,M}}~\sum\limits_{(x_p,y_p)\in \Omega_p}y_p},$$
where $Q(x_p,y_p)=p$ with $x_p>y_p$ and $Pr(N)$ denotes the $N$-th prime. 

\begin{figure}[H]
    \centering
    \includegraphics[width=10cm]{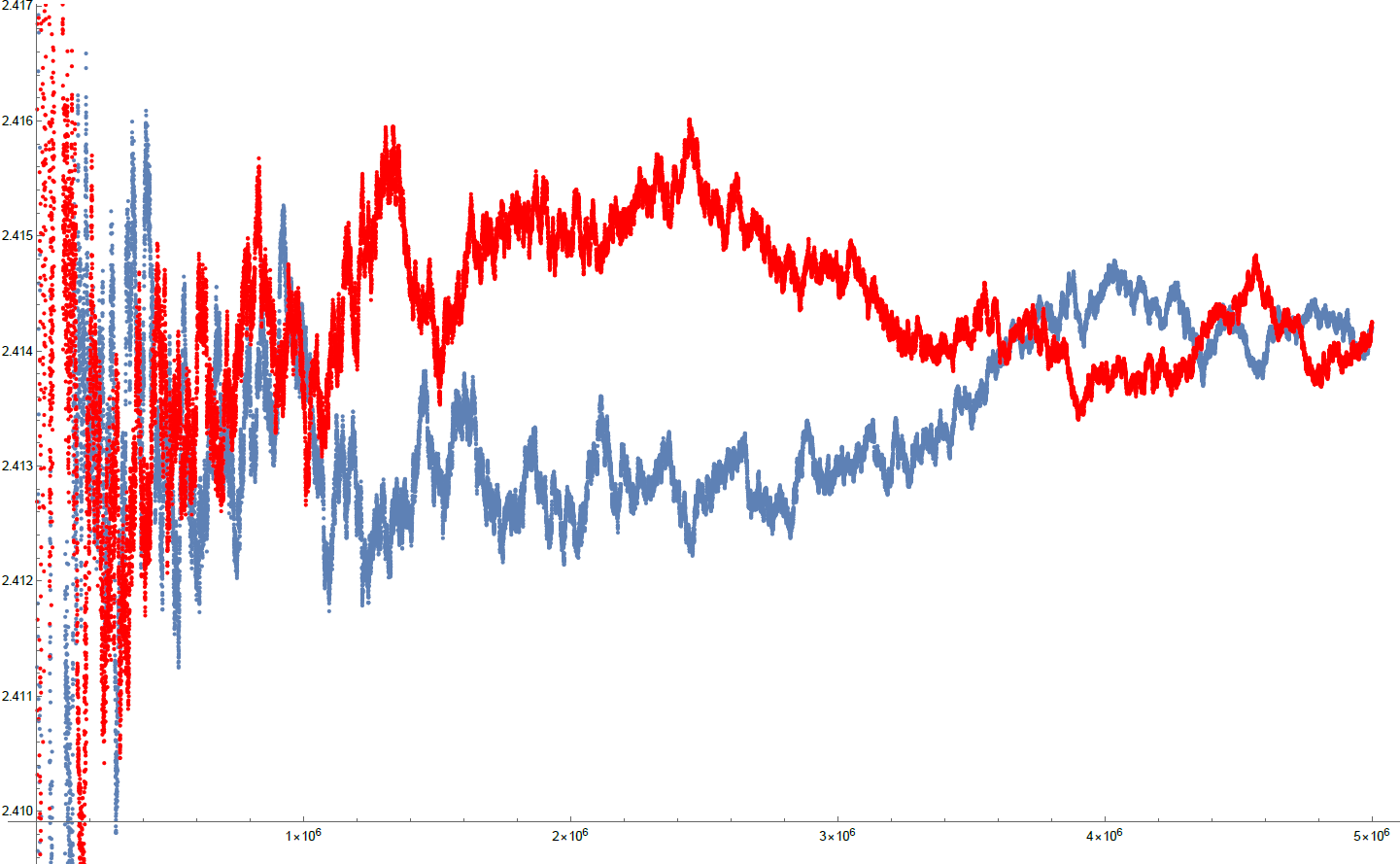}
    \caption{$\fr(N;8,1)$ and $\fr(N;8,5)$ for $N\in[100,5\times10^6]$:}
     $\fr(N;8,1)$ is marked blue, $\fr(N;8,5)$ is marked red. We only mark points with $100|N$.
    \label{p1-p5}
\end{figure}

For $Q(x,y)=x^2+y^2$, we have proved that $$\lim_{N\to\infty}\frac{\sum\limits_{p\leq N,p\equiv1\pmod{8}}x_p}{\sum\limits_{p\leq N,p\equiv1\pmod{8}}y_p}=\lim_{N\to\infty}\frac{\sum\limits_{p\leq N,p\equiv5\pmod{8}}x_p}{\sum\limits_{p\leq N,p\equiv5\pmod{8}}y_p}=1+\sqrt{2},$$
that is, $$\lim_{N\to\infty}\fr(N;8,1)=\lim_{N\to\infty}\fr(N;8,5)=1+\sqrt{2}.$$
To investigate the difference between  
$\fr(N;8,1)$ and $\fr(N;8,5)$, we compute these two functions with numerical data up to $5\times10^6$, and the graph is presented in Figure \ref{p1-p5}.

\begin{figure}[H]
    \centering
    \includegraphics[width=10cm]{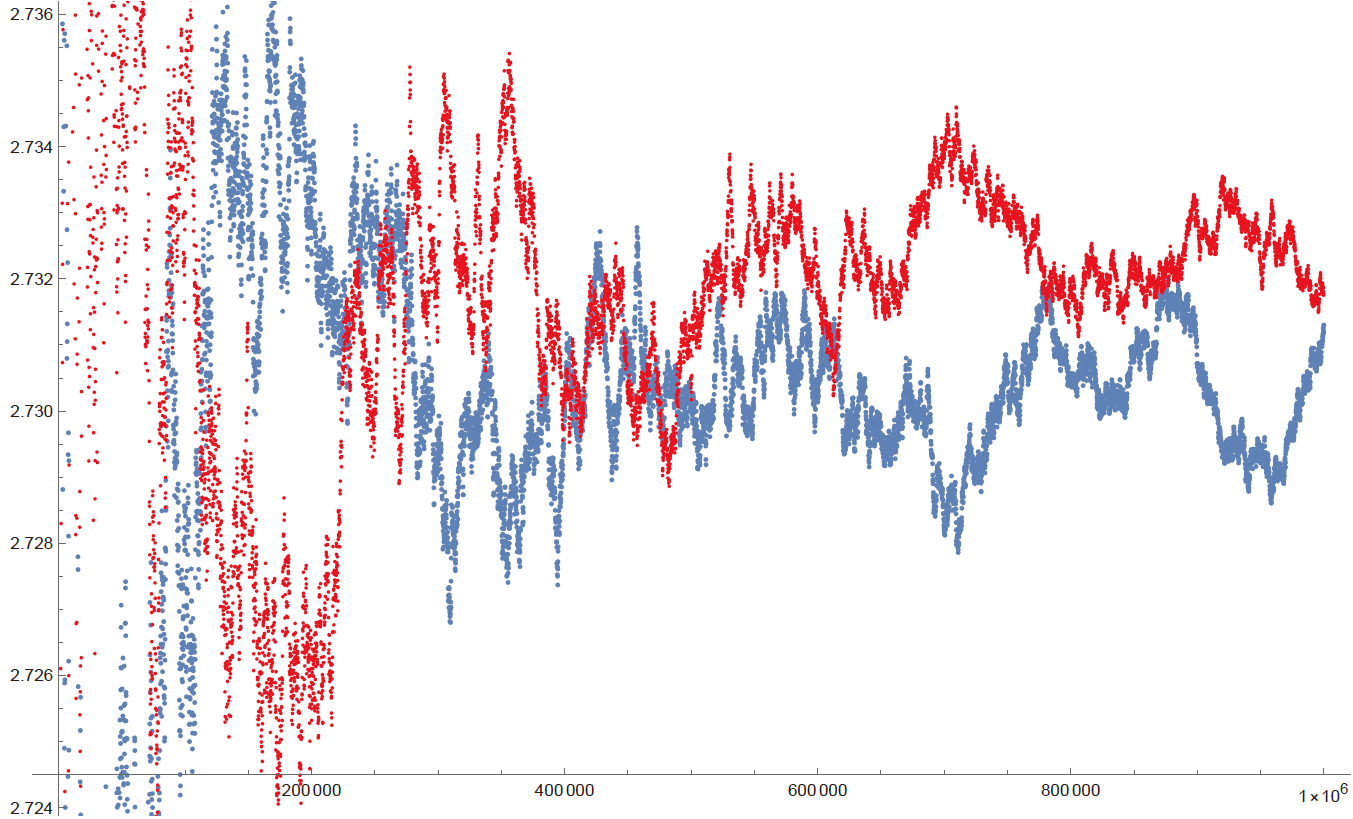}
    \caption{$\fr(N;12,1)$ and $\fr(N;12,7)$ for $N\in[100,1\times10^6]$:}
     $\fr(N;12,1)$ is marked blue, $\fr(N;12,7)$ is marked red.
    We only mark points with $100|N$. \label{1(12)vs7(12)}
\end{figure}

Similarly, for another quadratic form $Q(x,y)=x^2+xy+y^2$,
we consider $\fr(N;12,1)$ and $\fr(N;12,7)$ and compute these two functions with numerical data up to $1\times 10^6$.  Then we have the graph in Figure \ref{1(12)vs7(12)}.

We are interested in the phenomenon of oscillation and repeated intersections of two functions as are shown in the above two figures. Now we define another function, which generates a graph that preserves oscillatory behavior and intersections while exhibiting nice symmetry. This symmetry facilitates more effective observation and analysis. 

\begin{figure}[H]
    \centering
    \includegraphics[width=13cm]{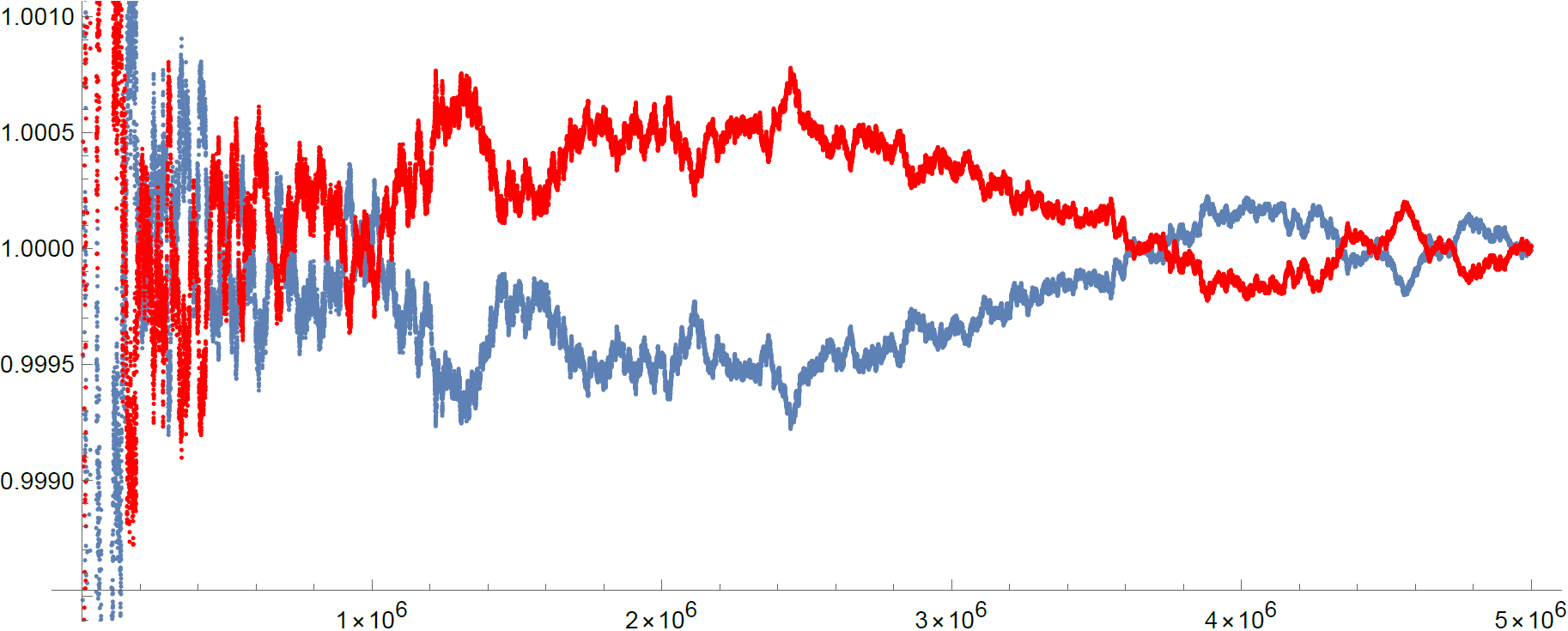}
    \caption{$R(N;8,1)$ and $R(N;8,5)$ for $N\in[100,5\times10^6]$:}
    $R(N;8,1)$ is marked blue, $R(N;8,5)$ is marked red. We only mark points with $100|N$. This graph shows very nice symmetry since $\fr(N;8,1)$, $\fr(N;8,5)$ and $\fr(N)$ have the same limit as $N\to\infty$. 
    \label{R1-R2}
\end{figure}

For a given quadratic form $Q(x,y)$, take $R(N;M,m)=\frac{\fr(N;M,m)}{\fr(N)}$, where $\fr(N)=\frac{\sum_{p\leq Pr(N)}~\sum_{(x_p,y_p)\in \Omega_p}x_p}{\sum_{p\leq Pr(N)}~\sum_{(x_p,y_p)\in \Omega_p}y_p}$. 
For $Q(x,y)=x^2+y^2$, we consider $R(N;8,1)$ and $R(N;8,5)$ and compute these two functions 
with numerical data up to $5\times 10^6$, and the graph is presented in Figure \ref{R1-R2}. 
These three figures above show that although the two functions converge to the same limit, they repeatedly intersect and alternately surpass each other as $N$ increases. Based on this phenomenon, we propose the following conjecture.

\begin{conjecture}
The two differences $\fr(N;8,1)-\fr(N;8,5) $ and $\fr(N;12,1)-\fr(N;12,7) $ will change signs infinitely many times as $N\to\infty$.

More generally, let $Q(x,y)=ax^2+bxy+cy^2$ be a primitive positive definite integral binary quadratic form with $b+c\neq0$. Let $m_1$ and $m_2$ be distinct reduced residue classes modulo $M$ satisfying the positive-density hypothesis of Theorem 1.3 for $Q$. Then the difference $F(N;M,m_1)-F(N;M,m_2)$ changes sign infinitely often as $N\to \infty$.
\end{conjecture}

\subsection{Counting functions}\label{4.2}
\begin{figure}[H]
    \centering
    \includegraphics[width=15cm]{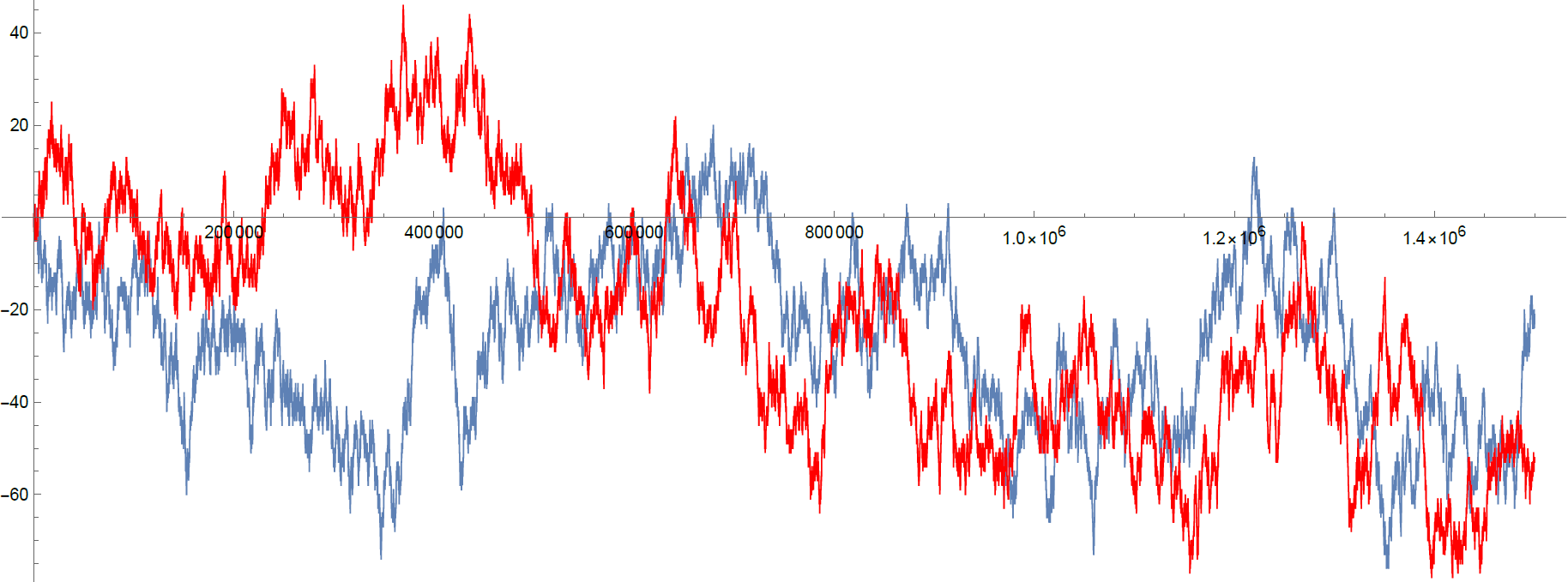}
    \caption{$D_1(x)$ and $D_5(x)$ for $x\in [2,1.5\times10^6]$:}
    $D_1(x)$ is marked blue, $D_5(x)$ is marked red.
    \label{D1D2}
\end{figure}

L. Devin \cite{devin2021discrepancies} conjectured that for more than half of $x>2$ in logarithmic scale, more than half of the primes below $x$ can be written as a sum of two squares with the even square larger than the odd one. Motivated by this conjecture, for $i=1,5$, we define 
\[
A_i(x)=\left\{p<x\mid p~\text{is prime},~p\equiv i\pmod8,~p=a^2+4b^2,~|a|>|2b|,~a,b\in \mathbb{Z}
\right\},\]
\[B_i(x)=\left\{p<x\mid p~\text{is prime},~p\equiv i\pmod8,~ p=a^2+4b^2,~|a|<|2b|,~a,b\in \mathbb{Z}\right\}.
\]
Recall that for a nonzero prime ideal $\mathfrak{p}$ of $\mathbb{Z}[\sqrt{-1}]$ prime to $2$, there exists a unique generator $\alpha$ of $\mathfrak{p}$ satisfying $\alpha \equiv 1 \mod(1+\sqrt{-1})^3$, and we define $\Lambda(\mathfrak{p})=\frac{\alpha}{|\alpha|}$ as in Remark \ref{remark2.4}. For each prime $p=a^2+4b^2$, there exist two prime ideals $\mathfrak{p}_1=(a+2b\sqrt{-1})$ and $\mathfrak{p}_2=(a-2b\sqrt{-1})$ lying over $p$, and exactly one of $\arg\left(\Lambda(\mathfrak{p}_1)\right)$ and $\arg\left(\Lambda(\mathfrak{p}_2)\right)$ lies in $(0,\frac{\pi}{2})\cup(\frac{\pi}{2},\pi)$. We correspond $p$ with this $\mathfrak{p}_i$ and let $\theta_p=\arg(\Lambda(\mathfrak{p}_i))$. Conversely, for a prime ideal $\mathfrak{p}$ prime to $2$ satisfying $\arg\Lambda(\mathfrak{p})\in (0,\frac{\pi}{2})\cup(\frac{\pi}{2},\pi)$, let $\alpha=a+b'\sqrt{-1}\equiv1\mod(1+\sqrt{-1})^3$ be the generator of $\mathfrak{p}$. Then $b'$ must be even and $\mathrm{N}(\mathfrak{p})=p=a^2+4b^2$, where $b'=2b$. This is a 1-1 correspondence, under which we have
\[
|a|>|2b|\quad\Longleftrightarrow\quad\theta_p\in(0,\frac{\pi}{4})\cup(\frac{3\pi}{4},\pi),\quad|a|<|2b|\quad\Longleftrightarrow\quad\theta_p\in (\frac{\pi}{4},\frac{\pi}{2})\cup(\frac{\pi}{2},\frac{3\pi}{4}).
\]
Remark \ref{remark2.4} implies that the angles \(\theta_p\) are equidistributed in $[0,\pi]$. Therefore, we have
\[
|A_i(x)|\sim|B_i(x)|\sim\frac{1}{2}\pi(x;8,i)\sim\frac{x}{8\log x},~i=1,5,
\]
where $\pi(x;8,i)=
\#\{p<x\mid p\text{ is prime and }p\equiv i\pmod 8\}$.
Define
\[
D_i(x)=|A_i(x)|-|B_i(x)|
\]
Then $D_i(x)=o\left(\frac{x}{\log x}\right)$ for $i=1,5$. Figure \ref{D1D2} shows how the two differences $D_1(x)$ and $D_5(x)$ vary with $x$. Theorem \ref{2.2}, together with Remark \ref{remark2.4}, does not determine whether this difference is positive or negative. To explore this question numerically, we compute $D_1(x)$ and $D_5(x)$ with numerical data up to $1.5\times10^6$, and the graph is presented in Figure \ref{D1D2}.

As shown in Figure \ref{D1D2}, the numerical values of $D_i(x)$ suggest a possible bias towards negative values. This observation leads to the following conjecture.

\begin{conjecture}
\label{conjecture_D1D2}
    There is a bias towards negative values in the distribution of the values of the functions $D_1(x)$ and $D_5(x)$.
\end{conjecture}

\section*{Acknowledgments}
The authors are supported by the National Natural Science Foundation of China (Nos. 12231009, 124B1010).

\end{document}